\documentclass[runningheads]{svmult}

\usepackage{amssymb}    %
\usepackage{epsfig}     %

\makeatletter
  
\def\ps@copyright{\let\@mkboth\@gobbletwo
  \def\@oddhead{}%
  \let\@evenhead\@oddhead
  \def\@oddfoot{\small\slshape
    \def\@tempa{0}
    \ifx\@volume\@tempa
         {}\/%
    \else
      Article published in \@jou@vol@pag\hfil\hbox{}\fi}%
  \let\@evenfoot\@oddfoot
}
  
\makeatother
  
\def\setboxz@h{\setbox\z@\hbox}
\def\wdz@{\wd\z@}
  
\def\binrel@#1{\begingroup
  \setboxz@h{\thinmuskip0mu
    \medmuskip\m@ne mu\thickmuskip\@ne mu
    \setbox\tw@\hbox{$#1\m@th$}\kern-\wd\tw@
    ${}#1{}\m@th$}%
  \edef\@tempa{\endgroup\let\noexpand\binrel@@
    \ifdim\wdz@<\z@ \mathbin
    \else\ifdim\wdz@>\z@ \mathrel
    \else \relax\fi\fi}%
  \@tempa
}
                                                                                
\DeclareRobustCommand{\boldsymbol}[1]{%
  \begingroup
  \let\@nomath\@gobble \mathversion{bold}%
  \math@atom{#1}{%
  \mathchoice%
    {\hbox{$\m@th\displaystyle#1$}}%
    {\hbox{$\m@th\textstyle#1$}}%
    {\hbox{$\m@th\scriptstyle#1$}}%
    {\hbox{$\m@th\scriptscriptstyle#1$}}}%
  \endgroup}
  
\def\math@atom#1#2{%
   \binrel@{#1}\binrel@@{#2}}
  
\def\Addots{\mathinner{%
        \mkern 1 mu \raise 1 pt \hbox{.}
        \mkern 2 mu \raise 4 pt \hbox{.}
        \mkern 2 mu \raise 7 pt \vbox{\kern 7 pt \hbox{.}}
        \mkern 1 mu}}
  
\makeatother

\newtheorem{algorithm}{Algorithm}

\newcommand{\beq}[2]{\begin{equation}
              \label{#1} 
                {#2}
           \end{equation}}
 
\newcommand{\Matrix}[1]{\left[ \matrix{#1} \right]}

\newcommand{\bea}{\begin{array}}
\newcommand{\ea}{\end{array}}
 
\def\q{\quad}
\def\mr#1{$(\ref{#1})$}
 
\newcommand{\ignore}[1]{} 
 
\newcommand{\real}{\mathbb{R}}
\newcommand{\complex}{\mathbb{C}}

\newcommand{\Span}{\mathop{\rm span}\nolimits}

\newcommand{\Range}{\mathop{\rm range}\nolimits}

\newcommand{\Ac}{\mathcal A}
\newcommand{\Bc}{\mathcal B}
\newcommand{\Cc}{\mathcal C}
\newcommand{\Dc}{\mathcal D}
\newcommand{\Ec}{\mathcal E}

\newcommand{\Jc}{\mathcal J}
\newcommand{\Kc}{\mathcal K}
\newcommand{\Lc}{\mathcal L}
\newcommand{\Mc}{\mathcal M}

\newcommand{\Oc}{\mathcal O}

\newcommand{\Rc}{\mathcal R}

\newcommand{\Vc}{\mathcal V}
\newcommand{\Wc}{\mathcal W}

\newcommand{\sinv}{\frac{1}{s}}

\begin{document}

\title*{Pad\'e-Type Model Reduction of Second-Order and
Higher-Order Linear Dynamical Systems}
\titlerunning{Pad\'e-Type Model Reduction of Higher-Order 
Linear Dynamical Systems}
\author{Roland W. Freund}%
\institute{Department of Mathematics,
   University of California at Davis,
   One Shields Avenue,
   Davis, CA 95616, U.S.A.\\
\texttt{freund@math.ucdavis.edu}}
\maketitle

\begin{abstract}
A standard approach to reduced-order modeling of higher-order 
linear dynamical systems is to rewrite the system as an
equivalent first-order system and then employ Krylov-subspace
techniques for reduced-order modeling of first-order systems.
While this approach results in reduced-order models that are
characterized as Pad\'e-type or even true Pad\'e approximants of
the system's transfer function, in general, these models do not preserve 
the form of the original higher-order system.
In this paper, we present a new approach to reduced-order modeling
of higher-order systems based on projections onto suitably
partitioned Krylov basis matrices that are obtained by applying
Krylov-subspace techniques to an equivalent first-order system.
We show that the resulting reduced-order models preserve the
form of the original higher-order system. 
While the resulting reduced-order models are no 
longer optimal in the Pad\'e sense, we show that they still satisfy a 
Pad\'e-type approximation property.  
We also introduce the notion of Hermitian higher-order linear
dynamical systems, and we establish an enhanced Pad\'e-type approximation 
property in the Hermitian case.
\end{abstract}

\section{Introduction} \label{sec-intro}

The problem of model reduction is to replace
a given mathematical model of a system or process by
a model that is much smaller than the original model,
yet still describes---at least approximately---certain
aspects of the system or process.
Model reduction involves a number of interesting 
issues.
First and foremost is the issue of selecting appropriate
approximation schemes that allow the definition of suitable
reduced-order models.
In addition, it is often important that the reduced-order
model preserves certain crucial properties of the original
system, such as stability or passivity.
Other issues include the characterization of the
quality of the models, the extraction of the data
from the original model that is needed to actually
generate the reduced-order models, and the efficient and
numerically stable computation of the models.
 
In recent years, there has been a lot of interest in 
model-reduction techniques based on Krylov subspaces; 
see, for example, the survey papers~\cite{Fre97b,Fre00a,Bai02,Fre03b}.
The development of these methods was motivated mainly by the
need for efficient reduction techniques in VLSI circuit
simulation.
An important problem in that application area is the reduction
of very large-scale RCL subcircuits that arise in the modeling
of the chip's wiring, the so-called \textit{interconnect}. 
In fact, many of the Krylov-subspace reduction techniques that
have been proposed in recent years are tailored to RCL subcircuits.

Krylov-subspace techniques can be applied directly only to
first-order linear dynamical systems.
However, there are important applications that lead to second-order,
or even general higher-order, linear dynamical systems.
For example, RCL subcircuits are actually second-order linear dynamical 
systems.
The standard approach to employing Krylov-subspace
techniques to the dimension reduction of a second-order or
higher-order system is to first rewrite the system as an
equivalent first-order system and then apply Krylov-subspace
techniques for reduced-order modeling of first-order systems.
While this approach results in reduced-order models that are
characterized as Pad\'e-type or even true Pad\'e approximants of
the system's transfer function, in general, these models do not preserve 
the form of the original higher-order system.

In this paper, we describe an approach to reduced-order modeling
of higher-order systems based on projections onto suitably
partitioned Krylov basis matrices that are obtained by applying
Krylov-subspace techniques to an equivalent first-order system.
We show that the resulting reduced-order models preserve the
form of the original higher-order system. 
While the resulting reduced-order models are no 
longer optimal in the Pad\'e sense, we show that they still satisfy a 
Pad\'e-type approximation property.  
We further establish an enhanced Pad\'e-type approximation 
property in the special case of Hermitian higher-order linear
dynamical systems.

The remainder of the paper is organized as follows.  
In Section~\ref{sec-2nd}, we review the formulations of general
RCL circuits as first-order and second-order linear dynamical
systems.   
In Section~\ref{sec-horder}, we present our general framework for
special second-order and higher-oder linear dynamical systems. 
In Section~\ref{sec-forder}, we consider the standard reformulation
of higher-order systems as equivalent first-order systems.      
In Section~\ref{sec-rom}, we discuss
some general concepts of dimension
reduction of special second-order and general higher-order systems
via dimension reduction of corresponding first-order systems.
In Section~\ref{sec-pade}, we review the concepts of block-Krylov subspaces
and Pad\'e-type reduced-order models.
In Section~\ref{sec-sym}, we introduce the notion of Hermitian higher-order linear
dynamical systems, and we establish an enhanced Pad\'e-type approximation 
property in the Hermitian case.
In Section~\ref{sec-sprim}, we present the SPRIM algorithm for
special second-order systems.
In Section~\ref{sec-examples}, we report results of 
some numerical experiments with the SPRIM algorithm.
Finally, in Section~\ref{sec-cremarks}, we mention some open problems 
and make some concluding remarks.

Throughout this paper the following notation is used.
Unless stated otherwise, all vectors and matrices are allowed to have real 
or complex entries.
For a complex number $\alpha$ or a complex matrix $M$,
we denote its complex conjugate by $\overline{\alpha}$
or $\overline{M}$, respectively.
For a matrix $M = \Matrix{m_{jk}}$, $M^T := \Matrix{m_{kj}}$
is the transpose of~$M$,
and $M^H :={\overline{M}}^T = \Matrix{\overline{m_{kj}}}$
is the conjugate transpose of~$M$.
For a square matrix $P$, we write $P \succeq 0$ if $P=P^H$ is
Hermitian and if $P$ is positive semidefinite, i.e.,
$x^H Px\ge 0$ for all vectors $x$ of suitable dimension.
We write $P \succ 0$ if $P=P^H$ is positive definite, i.e.,
$x^H Px > 0$ for all vectors $x$, except $x=0$.
The $n\times n$ identity matrix is denoted by $I_n$ and the
zero matrix by $0$. 
If the dimension of $I_n$ is apparent from the context,
we drop the index and simply use $I$.
The actual dimension of $0$ will always be
clear from the context.
The sets of real and complex numbers are denoted by $\real$ and
$\complex$, respectively.

\section{RCL circuits as first-order and second-order systems} \label{sec-2nd}

An important class of electronic circuits is linear RCL circuits that contain 
only resistors, capacitors, and inductors.
For example, such RCL circuits are used to model the interconnect of
VLSI circuits; see, e.g., \cite{CheLLC00,KimGP94,OdaCP98}.
In this section, we briefly review the RCL circuit equations and 
their formulations as first-order and second-order linear dynamical
systems.                                                                                

\subsection{RCL circuit equations}

General electronic circuits are usually modeled as networks whose branches
correspond to the circuit elements and whose nodes
correspond to the interconnections of the circuit elements;
see, e.g.,~\cite{VlaS94}.
Such networks are characterized by 
\textit{Kirchhoff's current law} (KCL),
\textit{Kirchhoff's voltage law} (KVL), and the 
\textit{branch constitutive relations} (BCRs).
The Kirchhoff laws depend only on the interconnections of the circuit elements,
while the BCRs characterize the actual elements.
For example, the BCR of a linear resistor is Ohm's law.
The BCRs are linear equations for simple devices, such as
linear resistors, capacitors, and inductors, and they
are nonlinear equations for more complex devices, such as
diodes and transistors.

The connectivity of such a network can be captured using a directional
graph.
More precisely, the nodes of the graph correspond to the nodes of the circuit,
and the edges of the graph correspond to each of the circuit elements.
An arbitrary direction is assigned to graph edges, so one can
distinguish between the source and destination nodes.  
The adjacency matrix, $A$, of the directional graph describes the connectivity 
of a circuit. Each row of $A$ corresponds to a graph edge and,
therefore, to a circuit element. 
Each column of $A$ corresponds to a graph or circuit node. 
The column corresponding to the datum (ground) node of the circuit is omitted in 
order to remove redundancy.
By convention, a row of $A$ contains $+1$ in
the column corresponding to the source node, $-1$ in the column
corresponding to the destination node, and $0$ everywhere else.
Kirchhoff's laws can be expressed in terms of $A$ as follows:
\beq{AAA}{
\bea{rrl}
\mbox{KCL:} & \q A^T i_b & = 0,   \\[4pt]
\mbox{KVL:} & \q   A v_n & = v_b.
\ea
}
Here, the vectors $i_b$ and $v_b$ contain the branch currents and
voltages, respectively, and $v_n$ the non-datum node
voltages.
 
We now restrict ourselves to linear RCL circuits, and for simplicity, we
assume that the circuit is excited only by current sources.
In this case, $A$, $v_b$, and $i_b$
can be partitioned according to circuit-element types as follows:
\beq{AAB}{
A   = \Matrix{ A_i \cr 
               A_g \cr 
               A_c \cr 
               A_l }, \q
v_b = v_b(t) = \Matrix{ v_i \cr 
                        v_g \cr
                        v_c \cr
                        v_l }, \q
i_b = i_b(t) = \Matrix{ i_i \cr 
                       i_g \cr
                       i_c \cr
                       i_l }.
}
Here, the subscripts $i$, $g$, $c$, and $l$ stand for branches
containing current sources, resistors, capacitors, and inductors,
respectively.
Using~\mr{AAB}, the KCL and KVL equations~\mr{AAA} take on the 
following form:
\beq{AAC}{
\bea{rl} 
A_i^T i_i + A_g^T i_g + A_c^T i_c + A_l^T i_l & = 0,\\[4pt]
A_i v_n = v_i,\q 
A_g v_n = v_g,\q 
A_c v_n = v_c,\q 
A_l v_n & = v_l.
\ea
}   
Furthermore, the BCRs can be stated as follows:
\beq{AAD}{
i_i = -I(t), \q i_g = G v_g, \q
i_c = C \frac{d}{dt} v_c, \q
v_l = L \frac{d}{dt} i_l.
}
Here, $I(t)$ is the vector of current-source values, $G \succ 0$ and
$C \succ 0$ are diagonal matrices whose
diagonal entries are the conductance and capacitance values of the
resistors and capacitors, respectively, and 
$L \succeq 0$ is the inductance matrix.
In the absence of inductive coupling, $L$ is also a diagonal
matrix, but in general, $L$ is a full matrix.
However, an important special case is inductance matrices $L$ whose
inverse, the so-called susceptance matrix,
$S = L^{-1}$ is sparse; see~\cite{ZheKBP02,ZheP02}. 

Equations~\mr{AAC} and~\mr{AAD}, together with initial conditions
for $v_n(t_0)$ and $i_l(t_0)$ at some initial time $t_0$, 
provide a complete description of a given RCL circuit.
For simplicity, in the following we assume $t_0=0$ with zero initial
conditions:
\beq{AAE}{
v_n(0) = 0\q \mbox{and}\q i_l(0) = 0.
}
Instead of solving~\mr{AAC} and~\mr{AAD} directly, one usually first eliminates 
as many variables as possible;
this procedure is called modified nodal analysis (MNA)~\cite{HoRB75,VlaS94}.
More precisely, using the last three equations in~\mr{AAC} and the
first three equations in~\mr{AAD}, one can eliminate $v_g$, $v_c$, $v_l$,
$i_i$, $i_g$, $i_c$, and is left with the coupled equations
\beq{AAF}{
\bea{rl}
A_i^T I(t) & = A_g^T G A_g v_n + A_c^T C A_c 
             \displaystyle{\frac{d}{dt}} v_n + A_l^T i_l,\\[6pt]
A_l v_n & = L \displaystyle{\frac{d}{dt}} i_l
\ea
}
for $v_n$ and $i_l$.
Note that the equations~\mr{AAF} are completed by the initial conditions~\mr{AAE}.

For later use, we remark that the energy supplied to the RCL circuit by the
current sources is given by
\beq{AAG}{
E(t) = \int_0^t \bigl(v_i(\tau)\bigr)^T I(\tau)\, d\tau.
}
Recall that the entries of the vector $v_i$ are the voltages at the current sources.
In view of the second equation in~\mr{AAC}, $v_i$ is connected to $v_n$ by
the output relation
\beq{AAH}{
v_i = A_i v_n.
}

\subsection{RCL circuits as first-order systems}

The RCL circuit equations~\mr{AAF} and~\mr{AAH} can be viewed as a 
first-order time-invariant linear dynamical system with state vector
\[
z(t) := \Matrix{v_n(t)\cr 
                i_l(t)},
\]
and input and output vectors
\beq{ABB}{
u(t) := I(t)\q \mbox{and} \q y(t) := v_i(t),
}
respectively.
Indeed, the equations~\mr{AAF} and~\mr{AAH} are equivalent to
\beq{ABC}{
\bea{rl} 
\Ec\, \displaystyle{\frac{d}{dt}} z(t) - \Ac\, z(t) & = \Bc\, u(t), \\[8pt]
y(t)  & = \Bc^T z(t),
\ea 
}
where
\beq{ABD}{
\Ec := \Matrix{A_c^T C A_c & 0 \cr
                        0 & L},\q
\Ac := \Matrix{-A_g^T G A_g & -A_l^T \cr
                       A_l & 0},\q
\Bc := \Matrix{A_i^T \cr 0}.
}
Note that~\mr{ABC} is a system of \textit{differential-algebraic equations} (DAEs) of
first order.
Furthermore, in view of~\mr{ABB}, the energy~\mr{AAG}, 
which is supplied to the RCL circuit by the current sources, is just
the integral
\beq{ABE}{
E(t) = \int_0^t \bigl(y(\tau)\bigr)^T u(\tau)\, d\tau
}
of the inner product of the input and output vectors of~\mr{ABC}.
RCL circuits are passive systems, i.e., they do not generate energy,
and~\mr{ABE} is an important formula for the proper treatment of passivity;
see, e.g.,~\cite{AndV73,LozBEM00}.

\subsection{RCL circuits as second-order systems}

In this subsection, we assume that the inductance matrix $L$ of the RCL
circuit is nonsingular.
In this case, the RCL circuit equations~\mr{AAF} and~\mr{AAH} can also be viewed as a 
second-order time-invariant linear dynamical system with state vector
\[ 
x(t) := v_n(t),
\]
and the same input and output vectors~\mr{ABB} as before.
Indeed, by integrating the second equation of~\mr{AAF} and using~\mr{AAE},
we obtain 
\beq{ACB}{
L\, i_l(t) = A_l \int_0^t v_n(\tau)\, d\tau.
}
Since $L$ is assumed to be nonsingular, we can employ the relation~\mr{ACB} 
to eliminate $i_l$ in~\mr{AAF}.
The resulting equation, combined with~\mr{AAH}, can be written
as follows:
\beq{ACC}{
\bea{rl} 
P_1 \, \displaystyle{\frac{d}{dt}} x(t) + P_0\, x(t) +
 P_{-1}\, \int_0^t x(\tau)\, d\tau & = B \, u(t), \\[8pt]
y(t)  & = B^T x(t).
\ea 
}
Here, we have set
\beq{ACD}{
P_1 := A_c^T C A_c ,\q
P_0 := A_g^T G A_g, \q
P_{-1} := A_l^T L^{-1} A_l,\q
B := A_i^T.
}
Note that the first equation in~\mr{ACC} is a system of integro-DAEs.
We will refer to~\mr{ACC} as a \textit{special} second-order time-invariant linear 
dynamical system.
We remark that the input and output vectors of~\mr{ACC} are the same as
in the first-order formulation~\mr{ABC}.
In particular, the important formula~\mr{ABE} for the energy supplied
to the system remains valid for the special second-order
formulation~\mr{ABC}.

If the input vector $u(t)$ is differentiable, then by differentiating the first
equation of~\mr{ACC} we obtain the ``true'' second-order formulation
\beq{ACF}{
\bea{rl} 
P_1 \, \displaystyle{\frac{d^2}{dt^2}} x(t) + P_0\, \displaystyle{\frac{d}{dt}} x(t) +
 P_{-1}\, x(t) & = B \, \displaystyle{\frac{d}{dt}} u(t), \\[8pt]
y(t)  & = B^T x(t).
\ea 
}
However, besides the additional assumption on the differentiability of $u(t)$, 
the formulation~\mr{ACF}
also has the disadvantage that the energy supplied to the system is no longer given
by the integral of the inner product of the input and output vectors
\[
\hat{u}(t) := \displaystyle{\frac{d}{dt}} u(t) \q \mbox{and}\q \hat{y}(t) := y(t)
\] 
of~\mr{ACF}.
Related to this lack of a formula of type~\mr{ABE} is the fact that the
transfer function of~\mr{ACF} is no longer positive real.
For these reasons, we prefer to use the special second-order formulation~\mr{ACC},
rather than the more common formulation~\mr{ACF}.

\section{Higher-order linear dynamical systems} \label{sec-horder}

In this section, we discuss our general framework for special second-order
and higher-order linear dynamical systems.
We denote by $m$ and $p$ the number of inputs and outputs,
respectively, and by $l$ the order of such systems.
In the following, the only assumption on $m$, $p$, and $l$ is
that $m$, $p$, $l \geq 1$.

\subsection{Special second-order systems} 

A \textit{special second-order $m$-input $p$-output time-invariant linear 
dynamical system of order $l$} is a system of integro-DAEs
of the following form:
\beq{ASA}{
\bea{rl}
P_1\, \displaystyle{\frac{d}{dt}} x(t)
    + P_0\, x(t) + P_{-1}\,
 \displaystyle{\int_{t_0}^t x(\tau)\, d\tau} & = B\, u(t), \\[8pt]
y(t)  & = D\, u(t) + L\, x(t), \\[8pt]
x(t_0) & = x_0.
\ea 
}
Here, $P_{-1}$, $P_0$, $P_1 \in \complex^{N \times N}$,
$B \in \complex^{N\times m}$, $D \in \complex^{p\times m}$,
and $L \in \complex^{p\times N}$ are given matrices,
$t_0 \in \real$ is a given initial time,
and $x_0 \in \complex^N$ is a given vector of initial values.
We assume that the $N\times N$ matrix
\[
s P_1 + P_0 + \sinv P_{-1}
\]
is singular only for finitely many values
of $s\in \complex$.

The frequency-domain transfer function of~\mr{ASA} is given by
\beq{ASG}{
H(s) = D + L \Bigl(s P_1 + P_0 + \sinv P_{-1} \Bigr)^{-1} B.
}
Note that
\[
H: \complex \mapsto \left(\complex \cup {\infty}\right)^{p \times m}
\]
is a matrix-valued rational function.

In practical applications, such as the case of RCL circuits described in
Section~\ref{sec-2nd}, the matrices $P_0$ and $P_1$ are
usually sparse.
The matrix $P_{-1}$, however, may be dense, but
has a sparse representation of the form
\beq{AF1}{
P_{-1} = F_1 G F_2^H
}
or
\beq{AF2}{
P_{-1} = F_1 G^{-1} F_2^H,\q \mbox{with  nonsingular $G$},
}
where $F_1$, $F_2\in \complex^{N \times N_0}$ and
$G\in \complex^{N_0 \times N_0}$ are
sparse matrices.
We stress that in the case~\mr{AF1}, the matrix $G$ is
not required to be nonsingular.
In particular, for any matrix $P_{-1} \in \complex^{N\times N}$,
there is always the trivial factorization~\mr{AF1} with
$F_1 = F_2 = I$ and $G = P_{-1}$.
Therefore, without loss of generality, in the following, we assume that the 
matrix $P_{-1}$ in~\mr{ASA} is given by a product of the form~\mr{AF1}
or~\mr{AF2}.
 
\subsection{General higher-order systems} 

An \textit{$m$-input $p$-output time-invariant linear dynamical
system of order $l$} is a system of DAEs of the following form:
\beq{AHA}{
\bea{rl}
& \quad\ P_l\, \displaystyle{\frac{d^l}{dt^l}} x(t)
 + P_{l-1}\, \displaystyle{\frac{d^{l-1}}{dt^{l-1}}} x(t)
 + \cdots + P_1\, \displaystyle{\frac{d}{dt}} x(t) + P_0\, x(t)
                   = B\, u(t), \\[8pt]
& y(t) = D\, u(t)
  + L_{l-1}\, \displaystyle{\frac{d^{l-1}}{dt^{l-1}}} x(t)
 + \cdots + L_1\, \displaystyle{\frac{d}{dt}} x(t) + L_0\, x(t).
\ea
}
Here, $P_i\in \complex^{N \times N}$, $0\leq i \leq l$,
$B \in \complex^{N\times m}$, $D \in \complex^{p\times m}$,
and $L_j \in \complex^{p\times N}$, $0\leq j < l$,
are given matrices, and $N$ is called the state-space dimension
of~\mr{AHA}.
Moreover, in~\mr{AHA},
 $u: [t_0, \infty) \mapsto \complex^m$
is a given input function, $t_0 \in \real$ is a given
initial time, the components of the vector-valued
function $x: [t_0, \infty) \mapsto \complex^N$ are the so-called
state variables, and $y: [t_0, \infty) \mapsto \complex^p$
is the output function.
The system is completed by initial conditions of the form
\beq{AHB}{
\displaystyle{\frac{d^i}{dt^i}} x(t) \biggm|_{t=t_0}
= x_0^{(i)},\q 0\leq i < l,
}
where $x_0^{(i)} \in \complex^N$, $0\leq i < l$, are given
vectors.
 
The frequency-domain transfer function of~\mr{AHA} is given
by
\beq{AHC}{
H(s) := D + L(s) \bigl(P(s)\bigr)^{-1} B,\q s \in \complex,
}
where
\beq{AHE}{
P(s) := s^l P_l + s^{l-1} P_{l-1} + \cdots + s P_1 + P_0
}
and
\[
L(s) := s^{l-1} L_{l-1} + s^{l-2} L_{l-2} + \cdots + s L_1 + L_0.
\]
Note that
\[
P: \complex \mapsto \complex^{N\times N}\q
\mbox{and}\q L: \complex \mapsto \complex^{p\times N}
\]
are matrix-valued polynomials,
and that
\[
H: \complex \mapsto \left(\complex \cup {\infty}\right)^{p \times m}
\]
is a matrix-valued rational function.
We assume that the polynomial~\mr{AHE}, $P$, is \textit{regular},
that is, the matrix $P(s)$ is singular only for finitely many values
of $s\in \complex$; see, e.g., \cite[Part II]{GohLR82}.
This guarantees that the transfer function~\mr{AHC} has only
finitely many poles.
 
\subsection{First-order systems}
 
For the special case $l=1$, systems of the form~\mr{AHA} are called
first-order systems.
In the following, we use calligraphic letters for the data
matrices and $z$ for the vector of state-space variables of
first-order systems.
More precisely, we always write first-order systems in the form
\beq{AFA}{
\bea{rl}
\Ec\, \displaystyle{\frac{d}{dt}} z(t) - \Ac\, z(t)
      & = \Bc\, u(t), \\[8pt]
y(t)  & = \Dc\, u(t) + \Lc\, z(t), \\[8pt]
z(t_0) & = z_0.  
\ea
}
Note that the transfer function of~\mr{AFA} is given by
\beq{AFB}{
H(s) = \Dc + \Lc\, \bigl(s\, \Ec - \Ac\bigr)^{-1} \Bc.
}

\section{Equivalent first-order systems} \label{sec-forder}

A standard approach to treat higher-order systems is to rewrite
them as equivalent first-order systems.
In this section, we present such equivalent first-order formulations
of special second-order and general higher-order systems.

\subsection{The case of special second-order systems}

We start with special second-order systems~\mr{ASA},
and we distinguish the two cases~\mr{AF1} and~\mr{AF2}.
 
First assume that $P_{-1}$ is given by~\mr{AF1}.
In this case, we set
\beq{CAG}{
   z_1(t) := x(t)\q \mbox{and}\q
   z_2(t) := F_2^H \displaystyle{\int_{t_0}^t x(\tau)\, d\tau}.
}
By~\mr{AF1} and~\mr{CAG}, the first relation in~\mr{ASA} can
be rewritten as follows:
\beq{CAI}{
P_1\, \displaystyle{\frac{d}{dt}} z_1(t)
  + P_0\, z_1(t) + F_1 G\, z_2(t) = B\, u(t).
}
Moreover,~\mr{CAG} implies that
\beq{CAK}{
G^H  \displaystyle{\frac{d}{dt}} z_2(t)
        = (F_2 G)^H z_1(t).
}
It follows from~\mr{CAG}--\mr{CAK} that the special second-order
system~\mr{ASA} (with $P_{-1}$ given by~\mr{AF1}) is equivalent
to a first-order system~\mr{AFA} where
\beq{CAM}{
\bea{rl}
z(t) & := \Matrix{z_1(t)\cr z_2(t)},\q
z_0 := \Matrix{x_0\cr 0},\q
\Lc := \Matrix{L & 0},\q
\Bc := \Matrix{B \cr 0},\\[15pt]
\Dc & := D,\q
\Ac := \Matrix{-P_0 & -F_1G \cr
               (F_2 G)^H & 0},\q
\Ec := \Matrix{P_1 & 0 \cr
                 0 & G^H}.
\ea
}
The state-space dimension of this first-order system is $N_1 := N + N_0$,
where $N$ and $N_0$ denote the dimensions of
$P_1 \in \complex^{N\times N}$ and $G \in \complex^{N_0 \times N_0}$.
Note that~\mr{AFB} is the corresponding representation of
the transfer function~\mr{ASG}, $H$, in terms of
the data matrices defined in~\mr{CAM}.
 
Next, we assume that $P_{-1}$ is given by~\mr{AF2}.
We set
\[
   z_1(t) := x(t)\q \mbox{and}\q
   z_2(t) := G^{-1} F_2^H \displaystyle{\int_{t_0}^t x(\tau)\, d\tau}.
\]
The first relation in~\mr{ASA} can then
be rewritten as
\[
P_1\,  \displaystyle{\frac{d}{dt}} z_1(t) + P_0\, z_1(t) + F_1\, z_2(t) = B\, u(t).
\]
Moreover, we have
\[
G\, \displaystyle{\frac{d}{dt}} z_2(t) = F_2^H z_1(t).
\]
It follows that the special second-order
system~\mr{ASA} (with $P_{-1}$ given by~\mr{AF2})
is equivalent to a first-order system~\mr{AFA} where
\beq{CAO}{ 
\bea{rl}
z(t) & := \Matrix{z_1(t)\cr z_2(t)},\q
z_0 := \Matrix{x_0\cr 0},\q
\Lc := \Matrix{L & 0},\q
\Bc := \Matrix{B \cr 0},\\[15pt]
\Dc & := D,\q
\Ac := \Matrix{-P_0 & -F_1 \cr
               F_2^H & 0},\q 
\Ec := \Matrix{P_1 & 0 \cr
                 0 & G}.
\ea
}
The state-space dimension of this first-order system is again $N_1 := N + N_0$.
Note that~\mr{AFB} is the corresponding representation of
the transfer function~\mr{ASG}, $H$, in terms of
the data matrices defined in~\mr{CAO}.

\subsection{The case of general higher-order systems}

It is well known (see, e.g., \cite[Chapter 7]{GohLR82}) that
any $l$-th order system with state-space dimension $N$ is equivalent
to a first-order system with state-space dimension~$N_1 := l N$.
Indeed, it is easy to verify that the $l$-th order system~\mr{AHA}
with initial conditions~\mr{AHB} is equivalent to the first-order 
system~\mr{AFA} with
\beq{CAC}{
\bea{rl}
&\!\!
        z(t) := \Matrix{x(t)\cr
                         \noalign{\vskip3pt}
                         \frac{d}{dt} x(t) \cr
                         \vdots \cr
                         \noalign{\vskip3pt}
                         \frac{d^{l-1}}{dt^{l-1}} x(t)},\q
z_0 := \Matrix{x_0^{(0)} \cr
              \noalign{\vskip3pt}
              x_0^{(1)} \cr
              \vdots \cr
              \noalign{\vskip3pt}
              x_0^{(l-1)}},\q
    \Bc := \Matrix{0 \cr \vdots \cr 0 \cr \noalign{\vskip3pt} B},\\[34pt]
&\!\!  
       \Lc := \Matrix{L_0 & L_1 & \cdots & L_{l-1}},\q \Dc := D,\\[12pt]
&\!\!
       \Ec := \Matrix{I  & 0 & 0 & \cdots & 0 \cr
               \noalign{\vskip2pt}
               0 & I & 0 & \cdots & 0\cr
              \vdots & \ddots & \ddots & \ddots & \vdots \cr
                  0 & \cdots & 0  & I  & 0 \cr
               \noalign{\vskip2pt}
                  0 & \cdots & 0 & 0 & P_l},\q
\Ac := -\Matrix{   0 & \!\! -I & 0 & \cdots & 0\cr
                  0 & 0 & \!\! -I & \ddots & \vdots \cr
             \vdots & \ddots &  \ddots  & \ddots  & 0 \cr
                  0 & \cdots & 0 & 0 & \!\! -I \cr
             \noalign{\vskip2pt}
 P_0 & P_1 & P_2 & \cdots & P_{l-1}}.
\ea
}
We remark that~\mr{AFB} is the corresponding representation of
the $l$-order transfer function~\mr{AHC}, $H$, in terms of
the data matrices defined in~\mr{CAC}.

\section{Dimension reduction of equivalent first-order systems}\label{sec-rom}

In this section, we discuss some general concepts of dimension
reduction of special second-order and general higher-order systems
via dimension reduction of equivalent first-order systems.

\subsection{General reduced-order models}

We start with general first-order systems~\mr{AFA}.
For simplicity, from now on we always assume zero initial conditions, i.e.,
$z_0 = 0$ in~\mr{AFA}.
We can then drop the initial conditions in~\mr{AFA}, and consider
first-order systems~\mr{AFA} of the following form:
\beq{AFA2}{
\bea{rl}
\Ec\, \displaystyle{\frac{d}{dt}} z(t) - \Ac\, z(t)
      & = \Bc\, u(t), \\[8pt]
y(t)  & = \Dc\, u(t) + \Lc\, z(t). \\[8pt]
\ea
}
Here, $\Ac$, $\Ec \in \complex^{N_1 \times N_1}$,
$\Bc_1 \in \complex^{N_1\times m}$, $\Dc \in \complex^{p\times m}$,
and $\Lc \in \complex^{p\times N_1}$ are given matrices.
Recall that $N_1$ is the state-space dimension of~\mr{AFA2}.
We assume that the matrix pencil $s\, \Ec - \Ac$ is \textit{regular},
i.e., the matrix  $s\, \Ec - \Ac$ is singular only for finitely many\
values of $s\in \complex$.
This guarantees that the transfer function of~\mr{AFA2},
\beq{DAB}{
H(s) := \Dc + \Lc\, \bigl(s\, \Ec - \Ac\bigr)^{-1} \Bc,
}
exists.
  
A \textit{reduced-order model} of~\mr{AFA2} is a system of the same
form as~\mr{AFA2}, but with smaller state-space dimension.
More precisely, a reduced-order model of~\mr{AFA2} with state-space
dimension $n_1$ ($<N_1$) is a system of the form 
\beq{DAC}{
\bea{rl}
\tilde{\Ec}\, \displaystyle{\frac{d}{dt}} \tilde{z}(t) - \tilde{\Ac}\, \tilde{z}(t)
      & = \tilde{\Bc}\, u(t), \\[8pt]
\tilde{y}(t)  & = \tilde{\Dc}\, u(t) +\tilde{\Lc}\, \tilde{z}(t),
\ea
}
where $\tilde{\Ac}$, $\tilde{\Ec} \in \complex^{n_1 \times n_1}$,
$\tilde{\Bc} \in \complex^{n_1\times m}$, $\tilde{\Dc} \in \complex^{p\times m}$,
and $\tilde{\Lc} \in \complex^{p\times n_1}$.
Again, we assume that the matrix pencil $s\, \tilde{\Ec} - \tilde{\Ac}$ is 
regular.
The transfer function of~\mr{DAC} is then given by
\beq{DAD}{
\tilde{H}(s) := \tilde{\Dc} + \tilde{\Lc}\,
       \bigl(s\, \tilde{\Ec} - \tilde{\Ac}\bigr)^{-1} \tilde{\Bc}.
}

Of course,~\mr{DAC} only provides a framework for model reduction.
The real problem, namely the choice of suitable matrices $\tilde{\Ac}$,
$\tilde{\Ec}$, $\tilde{\Bc}$, $\tilde{\Lc}$, $\tilde{\Dc}$, and 
sufficiently large reduced state-space dimension $n_1$ still remains to 
be addressed.

\subsection{Reduction via projection} \label{ssec-proj}

A simple, yet very powerful (when combined with Krylov-subspace
machinery) approach for constructing reduced-order models
is projection.
Let 
\beq{DAE}{
\Vc \in \complex^{N_1 \times n_1}
}
be a given matrix, and set
\beq{DAF}{
\tilde{\Ac} := \Vc^H \Ac\, \Vc,\q
\tilde{\Ec} := \Vc^H \Ec\, \Vc,\q
\tilde{\Bc} := \Vc^H \Bc\,\q
\tilde{\Lc} := \Lc\, \Vc,\q 
\tilde{\Dc} := \Dc.
}
Then, provided that the matrix pencil $s\, \tilde{\Ec} - \tilde{\Ac}$ 
is regular, the system~\mr{DAC} with matrices given by~\mr{DAF}
is a reduced-order model of~\mr{AFA2} with state-space dimension~$n_1$.

A more general approach employs two matrices,
\[
\Vc,\, \Wc \in \complex^{N_1 \times n_1},
\]
and two-sided projections of the form
\[
\tilde{\Ac} := \Wc^H \Ac\, \Vc,\q
\tilde{\Ec} := \Wc^H \Ec\, \Vc,\q
\tilde{\Bc} := \Wc^H \Bc\,\q
\tilde{\Lc} := \Lc\, \Vc,\q   
\tilde{\Dc} := \Dc.
\]
For example, the PVL algorithm~\cite{FelF94,FelF95} can be viewed as
a two-sided projection method, where the columns of the matrices
$\Vc$ and $\Wc$ are the first $n_1$ right and left Lanczos vectors
generated by the nonsymmetric Lanczos process~\cite{Lan50a}.

All model-reduction techniques discussed in the remainder of this
paper are based on projections of the type~\mr{DAF}.

Next, we discuss the application of projections~\mr{DAF} to 
first-order systems~\mr{AFA2} that arise as equivalent formulations
of special second-order and higher-oder linear dynamical
systems. 
Recall from Section~\ref{sec-forder} that such equivalent
first-order systems exhibit certain structures.
For general matrices~\mr{DAE}, $\Vc$, the projected matrices~\mr{DAF}
do not preserve these structures.
However, as we will show now, these structures are preserved
for certain types of matrices~$\Vc$. 

\subsection{Preserving special second-order structure} \label{ssec-pres}

In this subsection, we consider special second-order systems~\mr{ASA},
where $P_{-1}$ is either of the form~\mr{AF1} or~\mr{AF2}.
Recall that the data matrices of the equivalent first-order formulations
of~\mr{ASA} are defined in~\mr{CAM}, respectively~\mr{CAO}.

Let $\Vc$ be any matrix of the block form
\beq{DSA}{
\Vc = \Matrix{V_1 & 0 \cr
                     0  & V_2}, \q V_1 \in \complex^{N \times n},\q
                                   V_2 \in \complex^{N_0 \times n_0},\q
}
such that the matrix
\[
\tilde{G} := V_2^H G V_2\q \mbox{is nonsingular}.
\]
First, consider the case of matrices $P_{-1}$ of the form~\mr{AF1}.
Using~\mr{CAM} and~\mr{DSA}, one readily verifies that in this case,
the projected matrices~\mr{DAF} are as follows:
\beq{DSD}{
\bea{rl}
\tilde{\Ac} & = \Matrix{ -\tilde{P_0} & -\tilde{F}_1\tilde{G}\cr
  \bigr(\tilde{F}_2\tilde{G}\bigl)^H & 0},\q
\tilde{\Ec} = \Matrix{ \tilde{P_1} & 0\cr
                                 0  & \tilde{G}^H},\q
\tilde{\Bc} = \Matrix{\tilde{B} \cr 0},\\[15pt]
\tilde{\Lc} & = \Matrix{\tilde{L}  & 0},\q
\tilde{\Dc} = D.
\ea
}
Here, we have set
\beq{DSE}{
\tilde{P_0} := V_1^H P_0 V_1,\q
\tilde{P_1} := V_1^H P_1 V_1,\q
\tilde{B} :=  V_1^H B,\q 
\tilde{L} := L V_1,
}
and 
\[
\tilde{F}_1 := \bigl(V_1^H F_1 G V_2\bigr)\, \tilde{G}^{-1},\q
\tilde{F}_2 := \bigl(V_1^H F_2 G V_2\bigr)\, \tilde{G}^{-1}.
\]
Note that the matrices~\mr{DSD} are of the same form as
the matrices~\mr{CAM} of the first-order formulation~\mr{AFA2}
of the original special second-order system~\mr{ASA} (with $P_{-1}$
of the form~\mr{AF1}).
It follows that the matrices~\mr{DSD} define a
reduced-order model in special second-order form,
\beq{DSR}{
\bea{rl}
\tilde{P}_1\, \displaystyle{\frac{d}{dt}} \tilde{x}(t)
    + \tilde{P}_0\, \tilde{x}(t) + \tilde{P}_{-1}\,
 \displaystyle{\int_{t_0}^t \tilde{x}(\tau)\, d\tau} & =  \tilde{B}\, u(t), \\[8pt]
\tilde{y}(t)  & =  \tilde{D}\, u(t) +  \tilde{L}\, \tilde{x}(t), 
\ea
}
where
\[
\tilde{P}_{-1} := \tilde{F}_1 \tilde{G} \tilde{F}_2^H.
\]
We remark that the state-space dimension of~\mr{DSR} is $n$, where $n$ 
denotes the number of columns of the submatrix $V_1$ in~\mr{DSA}.

Next, consider the case of matrices $P_{-1}$ of the form~\mr{AF2}.
Using~\mr{CAO} and~\mr{DSA}, one readily verifies that in this case,
the projected matrices~\mr{DAF} are as follows:
\beq{DST}{
\bea{rl}
\tilde{\Ac} & = \Matrix{ -\tilde{P_0} & -\tilde{F}_1\cr
          \tilde{F}_2^H & 0},\q
\tilde{\Ec} = \Matrix{ \tilde{P_1} & 0\cr
                                 0  & \tilde{G}},\q
\tilde{\Bc} = \Matrix{\tilde{B} \cr 0},\\[15pt]
\tilde{\Lc} & = \Matrix{\tilde{L}  & 0},\q
\tilde{\Dc} = D.
\ea
}
Here, $\tilde{P_0}$, $\tilde{P_1}$, $\tilde{B}$, $\tilde{L}$ 
are the matrices defined in~\mr{DSE}, and
\[
\tilde{F}_1 := V_1^H F_1 V_2,\q
\tilde{F}_2 := V_1^H F_2 V_2.
\]
Again, the matrices~\mr{DST} are of the same form as
the matrices~\mr{CAO} of the first-order formulation~\mr{AFA2}
of the original special second-order system~\mr{ASA} (with $P_{-1}$
of the form~\mr{AF2}.
It follows that the matrices~\mr{DST} define a
reduced-order model in special second-order form~\mr{DSR}, where
\[
\tilde{P}_{-1} = \tilde{F}_1 \tilde{G}^{-1} \tilde{F}_2^H.
\]
 
\subsection{Preserving general higher-order structure} \label{ssec-preh}

We now turn to systems~\mr{AFA2} that are equivalent first-order 
formulations of general $l$-th order linear dynamical systems~\mr{AHA}.
More precisely, we assume that the matrices in~\mr{AFA2} are the ones
defined in~\mr{CAC}.

Let $\Vc$ be any $lN\times ln$ matrix of the block form  
\beq{SPC18}{
\Vc_n = \Matrix{S_n & 0 & 0 & \cdots & 0 \cr
                0 & S_n & 0 & \cdots & 0 \cr
                0 & 0 & \ddots & \ddots & \vdots \cr
           \vdots & \vdots & \ddots & \ddots & 0 \cr
                0 & 0 & \cdots & 0 & S_n},\q
S_n \in \complex^{N\times n},\q S_n^H S_n = I_n.
}
Although such matrices appear to be very special, they do arise in
connection with block-Krylov subspaces and lead to Pad\'e-type
reduced-order models; see Subsection \ref{ssec-pade} below.
The block structure~\mr{SPC18} implies that
the projected matrices~\mr{DAF} are given by
\beq{DCA}{ 
\bea{rl}
\tilde{\Ac} = -& \Matrix{   0 & \!\! -I & 0 & \cdots & 0\cr
                  0 & 0 & \!\! -I & \ddots & \vdots \cr
             \vdots & \ddots &  \ddots  & \ddots  & 0 \cr
                  0 & \cdots & 0 & 0 & \!\! -I \cr
             \noalign{\vskip2pt}
 \tilde{P}_0 & \tilde{P}_1 & \tilde{P}_2 & \cdots & \tilde{P}_{l-1}},\q 
       \tilde{\Ec} = \Matrix{I  & 0 & 0 & \cdots & 0 \cr
               \noalign{\vskip2pt}
               0 & I & 0 & \cdots & 0\cr
              \vdots & \ddots & \ddots & \ddots & \vdots \cr
                  0 & \cdots & 0  & I  & 0 \cr
               \noalign{\vskip2pt}
                  0 & \cdots & 0 & 0 & \tilde{P}_l},\\[42pt]
    \tilde{\Bc} = & \Matrix{0 \cr \vdots \cr 0 \cr 
             \noalign{\vskip3pt} \tilde{B}},\q
      \tilde{\Lc} = \Matrix{\tilde{L}_0 & \tilde{L}_1 & \cdots & \tilde{L}_{l-1}},\q
      \tilde{\Dc} = \Dc,
\ea
}
where 
\[
\tilde{P}_i := S_n^H P_i S_n,\q 0\leq i \leq l,\q
\tilde{B} := S_n^H B,\q
\tilde{L}_j := L_j S_n,\q 0\leq j < l.
\]
It follows that the matrices~\mr{DCA} define a
reduced-order model in $l$-th order form,
\beq{DCD}{
\bea{rl}
& \ \tilde{P}_l\, \displaystyle{\frac{d^l}{dt^l}} \tilde{x}(t)
  + \tilde{P}_{l-1}\, \displaystyle{\frac{d^{l-1}}{dt^{l-1}}} \tilde{x}(t)
  + \cdots + \tilde{P}_1\, \displaystyle{\frac{d}{dt}} \tilde{x}(t) 
  + \tilde{P}_0\, \tilde{x}(t) = \tilde{B}\, u(t), \\[8pt]
& \tilde{y}(t) = \tilde{D}\, u(t)
  + \tilde{L}_{l-1}\, \displaystyle{\frac{d^{l-1}}{dt^{l-1}}} \tilde{x}(t)
  + \cdots + \tilde{L}_1\, \displaystyle{\frac{d}{dt}} \tilde{x}(t) 
  + \tilde{L}_0\,  \tilde{x}(t),
\ea
}
of the original $l$-th order system~\mr{AHA}.
We remark that the state-space dimension of~\mr{DCD} is $n$, where $n$ 
denotes the number of columns of the matrix $S_n$ in~\mr{SPC18}.
 
\section{Block-Krylov subspaces and Pad\'e-type models} \label{sec-pade}

In this section, we review the concepts of block-Krylov subspaces
and Pad\'e-type reduced-order models.

\subsection{Pad\'e-type models}

Let $s_0 \in \complex$ be any point such that the matrix 
$s_0\, \Ec - \Ac$ is nonsingular.
Recall that the matrix pencil $s\, \Ec - \Ac$ is assumed to be
regular, and so the matrix $s_0\, \Ec - \Ac$ is nonsingular except
for finitely many values of $s_0 \in \complex$.
In practice, $s_0 \in \complex$ is chosen such that $s_0\, \Ec - \Ac$ is nonsingular
and at the same time, $s_0$ is in some sense ``close'' to a problem-specific
relevant frequency range in the complex Laplace domain.
Furthermore, for systems with real matrices $\Ac$ and $\Ec$
one usually selects $s_0 \in \real$ in order to avoid complex arithmetic.

We consider first-order systems of the form~\mr{AFA2} and their reduced-order
models of the form~\mr{DAC}.
By expanding the transfer function~\mr{DAB}, $H$, of the original system~\mr{AFA2}
about $s_0$, we obtain
\beq{PTR1}{
\bea{rl}
H(s) = \Lc\, \bigl( s\, \Ec - \Ac \bigr)^{-1} \Bc
     & = \Lc\, \Bigl( I + (s - s_0) \Mc \Bigr)^{-1} \Rc\\[4pt]
     & = \displaystyle{\sum_{i=0}^{\infty} (-1)^i \Lc\, \Mc^i\, \Rc\,
            (s-s_0)^i},
\ea
}
where
\beq{PTR2}{
\Mc := \bigl( s_0\, \Ec - \Ac \bigr)^{-1} \Ec
\q \mbox{and}\q
\Rc := \bigl( s_0\, \Ec - \Ac \bigr)^{-1} \Bc.
}
Provided that the matrix $s_0\, \tilde{\Ec} - \tilde{\Ac}$ is nonsingular,
we can also expand the transfer function~\mr{DAD}, $\tilde{H}$,
of the reduced-order model~\mr{DAC} about $s_0$.
This gives
\beq{PTR3}{
\bea{rl}
\tilde{H}(s) & = 
  \tilde{\Lc}\, \bigl( s\,  \tilde{\Ec} -  \tilde{\Ac} \bigr)^{-1} \Bc\\[4pt]
     & = \displaystyle{\sum_{i=0}^{\infty} (-1)^i \tilde{\Lc}\, \tilde{\Mc}^i\, 
           \tilde{\Rc}\, (s-s_0)^i}, 
\ea
} 
where
\[
\tilde{\Mc} := \bigl( s_0\, \tilde{\Ec} - \tilde{\Ac} \bigr)^{-1} \tilde{\Ec}
\q \mbox{and}\q
\Rc := \bigl( s_0\, \tilde{\Ec} - \tilde{\Ac} \bigr)^{-1} \tilde{\Bc}.
\]
We call the reduced-order model~\mr{DAC} a \textit{Pad\'e-type model} (with
expansion point $s_0$) of the original system~\mr{AFA2} if the 
Taylor expansions~\mr{PTR1} and~\mr{PTR3} agree in a number of leading terms, i.e.,
\beq{PTR5}{
\tilde{H}(s) = H(s) + {\Oc}\bigl((s-s_0)^{q} \bigr)  
}
for some $q = q(\tilde{\Ac},\tilde{\Ec},\tilde{\Bc},\tilde{\Lc},\tilde{\Dc}) > 0$.

Recall that the state-space dimension of the reduced-order model~\mr{DAC} is $n_1$.
If for a given $n_1$, the matrices $\tilde{\Ac}$, 
$\tilde{\Ec}$, $\tilde{\Bc}$, $\tilde{\Lc}$, $\tilde{\Dc}$ in~\mr{DAC} are chosen
such that $q=q(n_1)$ in~\mr{PTR5} is optimal, i.e., as large as possible, then  the 
reduced-order model~\mr{DAC} is called a \textit{Pad\'e model}.
All the reduced-order models discussed in the remainder of this paper are
Pad\'e-type models, but they are not optimal in the Pad\'e sense.

The (matrix-valued) coefficients in the expansions~\mr{PTR1} and~\mr{PTR3} are often 
referred to as \textit{moments}.
Strictly speaking, the term ``moments'' should only be used when $s_0 = 0$;
in this case, the Taylor coefficients of the Laplace-domain
transfer function directly correspond to the moments in time domain.
However, the use of the term ``moments'' has become common even in the
case of general $s_0$.
Correspondingly, the property~\mr{PTR5} is now generally referred to as
``moment matching''.

We remark that the moment-matching property~\mr{PTR5} is important for the
following two reasons.
First, for large-scale systems, the matrices $\Ac$ and $\Ec$ are usually
sparse, and the dominant computational work for moment-matching reduction
techniques is the computation of a sparse LU factorization of the
matrix $s_0\, \Ec - \Ac$.
Note that such a factorization is required already even if one only
wants to evaluate the transfer function $H$ at the point $s_0$.
Once a sparse LU factorization of $s_0\, \Ec - \Ac$ has been generated,
moments can be computed cheaply.
Indeed, in view of~\mr{PTR1} and~\mr{PTR2}, only sparse back solves, sparse matrix 
products
(with $\Ec$), and vector operations are required.
Second, the moment-matching property~\mr{PTR5} is inherently connected
to block-Krylov subspaces.
In particular, Pad\'e-type reduced-order models can be computed
easily be combining Krylov-subspace machinery and projection techniques.
In the remainder of the section, we describe this connection with 
block-Krylov subspaces.

\subsection{Block-Krylov subspaces}

In this subsection, we review the concept of block-Krylov
subspaces induced by the matrices $\Mc$ and $\Rc$ defined in~\mr{PTR2}.
Recall that $\Ac$, $\Ec \in \complex^{N_1 \times N_1}$ 
and $\Bc \in \complex^{N_1 \times m}$.
Thus we have
\beq{KKX}{
\Mc \in \complex^{N_1 \times N_1}\q \mbox{and}\q
\Rc \in \complex^{N_1 \times m}.
}
Next, consider the infinite \textit{block-Krylov matrix}
\beq{KKA}{ 
\Matrix{\Rc & \Mc\, \Rc & \Mc^2\, \Rc & \cdots & \Mc^j\, \Rc & \dots\, }.
}
In view of~\mr{KKX}, the columns of the matrix~\mr{KKA} are vectors 
in $\complex^{N_1}$,
and so only at most $N_1$ of these vectors are linearly independent.
By scanning the columns of the matrix~\mr{KKA} from left to right and
deleting each column that is linearly dependent on columns to its left,
one obtains the so-called \textit{deflated} finite block-Krylov matrix
\beq{KKC}{
\Matrix{\Rc^{(1)} & \Mc\, \Rc^{(2)} & \Mc^2\, \Rc^{(3)} & \cdots &
              \Mc^{j_{\max}-1}\, \Rc^{(j_{\max})}},
}
where each block $\Rc^{(j)}$ is a subblock of  $\Rc^{(j-1)}$, $j=1,2,\dots,j_{\max}$,
and $\Rc^{(0)} := \Rc$.
Let $m_j$ denote the number of columns of the $j$-th block $\Rc^{(j)}$.
Note that by construction, the matrix~\mr{KKC} has full column rank.
The \textit{$n$-th block-Krylov subspace} (induced by $\Mc$ and $\Rc$)
$\Kc_n\bigl(\Mc,\Rc\bigr)$ is defined as the subspace of $\complex^{N_1}$ spanned 
by the first $n$ columns of the matrix~\mr{KKC}; see,~\cite{AliBFH00} for more details 
of this construction.

Here, we will only use those block-Krylov subspaces that correspond to the
end of the blocks in~\mr{KKC}.
More precisely, let $n$ be of the form
\beq{KKE}{
n = n(j) := m_1 + m_2 + \cdots + m_j, \q \mbox{where}\q 1 \leq j \leq j_{\max}.
}
In view of the above construction, the $n$-th block-Krylov subspace is given by
\beq{KKG}{
\Kc_n\bigl(\Mc,\Rc\bigr) = 
  \Range\Matrix{\Rc^{(1)} & \Mc\, \Rc^{(2)} & 
                \Mc^2\, \Rc^{(3)} & \cdots & \Mc^{j-1}\, \Rc^{(j)}}.
}

\subsection{The projection theorem revisited} \label{ssec-prot}

It is well known that the projection approach described in Subsection~\ref{ssec-proj}
generates Pad\'e-type reduced-order models, provided that the matrix $\Vc$ in~\mr{DAE}
is chosen as a basis matrix for the block-Krylov subspaces induced by the
matrices $\Mc$ and $\Rc$ defined in~\mr{PTR2}.
This result is called the projection theorem, and it goes back to at 
least~\cite{VilS87}.
It was also established in \cite{Oda96,OdaCP97,OdaCP98} in connection with
the PRIMA reduction approach; see \cite{Fre00a} for more details.

One key insight to obtain structure-preserving Pad\'e-type reduced-order
models via block-Krylov subspaces and projection is the fact that the 
projection theorem remains valid when the above assumption on $\Vc$ is
replaced by the weaker condition
\beq{PTR9}{
\Kc_n\bigl(\Mc,\Rc\bigr) \subseteq \Range \Vc_n .
}
In this subsection, we present an extension of the projection theorem (as
stated in \cite{Fre00a}) to the case~\mr{PTR9}.

From now on, we always assume that $n$ is an integer of the form~\mr{KKE}
and that 
\beq{XXX9}{
\Vc_n \in \complex^{N_1 \times n_1} 
}
is a matrix satisfying~\mr{PTR9}. 
Note that~\mr{PTR9} implies $n_1 \geq n$.
We stress that we make no further assumptions about $n_1$.
We consider projected models given by~\mr{DAF} with $\Vc = \Vc_n$.
In order to indicate the dependence on the dimension $n$ of the
block-Krylov subspace in~\mr{PTR9}, we use the notation
\beq{XXX1}{
\bea{rl}
\Ac_n & := \Vc_n^H \Ac\, \Vc_n,\q 
\Ec_n := \Vc_n^H \Ec\, \Vc_n,\q
\Bc_n := \Vc_n^H \Bc,\\[6pt]
\Lc_n & := \Lc\, \Vc_n,\q 
\Dc_n := \Dc
\ea
}
for the matrices defining the projected reduced-order model.
In addition to~\mr{PTR9}, we also assume that
the matrix pencil $s\, \Ec_n - \Ac_n$ is regular,
and that at the expansion point $s_0$, the matrix $s_0\, \Ec_n - \Ac_n$ is 
nonsingular.
Then the reduced-order transfer function
\beq{XXX2}{
\bea{rl}
H_n(s) & := \Lc_n\, \bigl( s\, \Ec_n - \Ac_n \bigr)^{-1} \Bc_n\\[4pt]
     & = \Lc_n\, \Bigl( I + (s - s_0) \Mc_n \Bigr)^{-1} \Rc_n\\[4pt]
     & = \displaystyle{\sum_{i=0}^{\infty} (-1)^i 
       \Lc_n\, \Mc_n^i\, \Rc_n\, (s-s_0)^i}
\ea 
}
is a well-defined rational function.
Here, we have set
\beq{PTR8}{
\Mc_n := \bigl( s_0\, \Ec_n - \Ac_n \bigr)^{-1} \Ec_n
\q \mbox{and}\q
\Rc_n := \bigl( s_0\, \Ec_n - \Ac_n \bigr)^{-1} \Bc_n.
}
We remark that the regularity of the matrix pencil $s\, \Ec_n - \Ac_n$ implies that
the matrix $\Vc_n$ must have full column rank.

After these preliminaries, the extension of the 
projection theorem can be stated as follows.

\begin{theorem} \label{thm1}
Let\/~$n = n(j)$ be of the form\/~\mr{KKE}, and let\/~$\Vc_n$ be
a matrix satisfying\/~\mr{PTR9}.
Then the reduced-order model defined by\/~\mr{XXX1}
is a Pad\'e-type model with
\beq{XXX3}{
H_n(s) = H(s) + {\Oc}\bigl((s-s_0)^{j} \bigr). 
}
\end{theorem}
 
\begin{proof}
In view of~\mr{PTR1} and~\mr{XXX2}, the claim~\mr{XXX3} holds true if
\beq{PTR10}{
\Mc^i\, \Rc = \Vc_n\, \Mc_n^i\, \Rc_n\q \mbox{for all}\q i=0,1,\dots,j-1,
}
and thus it is sufficient to show~\mr{PTR10}.

By~\mr{KKG} and ~\mr{PTR9}, for each $i=0,1,\dots,j-1$, there exists
a matrix $\rho_i$ such that
\beq{PTR11}{ 
\Mc^{i}\, \Rc = \Vc_n\, \rho_i.
}
Moreover, since $\Vc_n$ has full column rank, each matrix $\rho_i$
is unique.
In fact, we will show that the matrices $\rho_i$
in~\mr{PTR11} are given by
\beq{PTR12}{
\rho_i = \Mc_n^i\, \Rc_n,\q i=0,1,\dots,j-1.
} 
The claim~\mr{PTR10} then follows by inserting~\mr{PTR12} 
into~\mr{PTR11}.

We prove~\mr{PTR12} by induction on $i$.
Let $i=0$.
In view of~\mr{PTR2} and~\mr{PTR11}, we have
\beq{PTR13}{
\Vc_n\, \rho_0 = \Rc = \bigl( s_0\, \Ec - \Ac \bigr)^{-1} \Bc.
}
Multiplying~\mr{PTR13} from the left by 
\beq{PTR15}{
\bigl( s_0\, \Ec_n - \Ac_n \bigr)^{-1} 
  \Vc_n^H \bigl( s_0\, \Ec - \Ac \bigr)
}
and using the definition of $\Rc_n$ in~\mr{PTR8}, it follows
that $\rho_0 = \Rc_n$.
This is just the relation~\mr{PTR12} for $i=0$.

Now let $1\leq i \leq j-1$, and assume that
\beq{PTR16}{
\rho_{i-1} = \Mc_n^{i-1}\, \Rc_n.
}
Recall that $\rho_{i-1}$ satisfies the equation~\mr{PTR11} (with $i$ 
replaced by $i-1$), and thus we have 
$\Mc^{i-1}\, \Rc = \Vc_n\, \rho_{i-1}$.
Together with~\mr{PTR11} and~\mr{PTR16}, it follows that
\beq{PTR17}{
\Vc_n\, \rho_i = \Mc^{i}\, \Rc = \Mc\, \bigl(\Mc^{i-1}\, \Rc \bigr)
 = \Mc\, \bigl(\Vc_n\, \rho_{i-1} \bigr)
 = \Mc\, \Vc_n \bigl(\Mc_n^{i-1}\, \Rc_n \bigr).
}
Note that, in view of the definition of $\Mc$ in~\mr{PTR2}, we have
\beq{PTR18}{
\Vc_n^H \bigl( s_0\, \Ec - \Ac \bigr) \Mc\, \Vc_n = 
\Vc_n^H \Ec\, \Vc_n = \Ec_n.
}
Multiplying~\mr{PTR17} from the left by the matrix~\mr{PTR15}
and using~\mr{PTR18} as well as the definition of $\Mc_n$ 
in~\mr{PTR8}, we obtain 
\[
\rho_i 
= \bigl( s_0\, \Ec_n - \Ac_n \bigr)^{-1} \Ec_n \,
  \bigl(\Mc_n^{i-1}\, \Rc_n \bigr)
= \Mc_n^{i}\, \Rc_n. 
\]
Thus the proof is complete.
\end{proof}

\subsection{Structure-preserving Pad\'e-type models} \label{ssec-pade}

We now turn to structure-preserving Pad\'e-type models.
Recall that, in Subsections~\ref{ssec-pres} and~\ref{ssec-preh}, we have
shown how special second-order and general higher-order structure is
preserved by choosing projection matrices of the form~\mr{DSA} 
and~\mr{SPC18}, respectively.
Moreover, in Subsection~\ref{ssec-prot} we pointed out that 
projected models are Pad\'e-type models if~\mr{PTR9} is satisfied.
It follows that the reduced-order models given by the projected data 
matrices~\mr{XXX1} are structure-preserving Pad\'e-type models, provided 
that the matrix $\Vc_n$ in \mr{XXX9}
is of the form~\mr{DSA}, respectively~\mr{SPC18}, and at the same time
fulfills the condition~\mr{PTR9}.  Next we show how to construct such
matrices $\Vc_n$.

Let 
\beq{QQQ1}{
\hat{\Vc}_n \in \complex^{N_1 \times n}
}
be any matrix whose columns span the $n$-th block-Krylov subspace 
$\Kc_n\bigl(\Mc,\Rc\bigr)$, i.e.,
\beq{QQQ2}{
\Kc_n\bigl(\Mc,\Rc\bigr) = \Range \hat{\Vc}_n.
}
First, consider the case of special second-order systems~\mr{ASA},
where $P_{-1}$ is either of the form~\mr{AF1} or~\mr{AF2}.
In this case, we partition $\hat{\Vc}_n$ as follows:
\beq{QQQ3}{ 
\hat{\Vc}_n = \Matrix{V_1\cr V_2},\q
 V_1 \in \complex^{N \times n},\q
                                   V_2 \in \complex^{N_0 \times n}.
}
Using the blocks in~\mr{QQQ3}, we set
\beq{QQQ4}{
\Vc_n := \Matrix{V_1 & 0 \cr 
                     0  & V_2}.
}
Clearly, the matrix~\mr{QQQ4} is of the form~\mr{DSA}, and
thus the projected models generated with $\Vc_n$ preserve
the special second-order structure.
Moreover, from~\mr{QQQ2}--\mr{QQQ4}, it follows that
\[
\Kc_n\bigl(\Mc,\Rc\bigr) = \Range \hat{\Vc}_n \subseteq \Range \Vc_n, 
\]
and so condition~\mr{PTR9} is satisfied.
Thus, the projected models are Pad\'e-type models and preserve second-order
structure.

Next, we turn to the case of general higher-order systems~\mr{AHA}.
In \cite{Fre04a}, we have shown that in this case, the block-Krylov subspaces induced
by the matrices $\Mc$ and $\Rc$, which are given by~\mr{CAC} and \mr{PTR2},
exhibit a very special structure. 
More precisely, the $n$-dimensional subspace $\Kc_n\bigl(\Mc,\Rc\bigr)$ 
of $\complex^{l N}$ can be viewed as $l$ copies of an  $n$-dimensional subspace
of $\complex^{N}$.
Let $S_n \in \complex^{N \times n}$ be a matrix whose columns form an orthonormal
basis of this  $n$-dimensional subspace
of $\complex^{N}$, and set
\beq{QQQ6}{
\Vc_n := \Matrix{S_n & 0 & 0 & \cdots & 0 \cr
                0 & S_n & 0 & \cdots & 0 \cr
                0 & 0 & \ddots & \ddots & \vdots \cr
           \vdots & \vdots & \ddots & \ddots & 0 \cr 
                0 & 0 & \cdots & 0 & S_n}.
}
From the above structure of the $n$-dimensional subspace $\Kc_n\bigl(\Mc,\Rc\bigr)$,
it follows that $\Vc_n$ satisfies the condition~\mr{PTR9}.
Furthermore, the matrix $\Vc_n$ is of the form~\mr{SPC18}.
Thus, the projected models generated with $\Vc_n$ are Pad\'e-type models 
and preserve higher-order structure.

In the remainder of this paper, we assume that $\Vc_n$ are matrices given
by~\mr{QQQ4} in the case of special second-order systems, 
respectively~\mr{QQQ6} in the case of higher-order systems, and we
consider the corresponding structure-preserving reduced-order models with
data matrices given by~\mr{XXX1}.

\section{Higher accuracy in the Hermitian case} \label{sec-sym}

For the structure-preserving Pad\'e-type models introduced in Subsection~\ref{ssec-pade},
the result of Theorem~\ref{thm1} can be improved further, provided the
underlying special second-order or higher-order linear dynamical system
is Hermitian, and the expansion point $s_0$ is real, i.e.,
\beq{HAA}{
s_0 \in \real.
}
More precisely, in the Hermitian case, the Pad\'e-type models obtained
via $\Vc_n$ match $2j(n)$ moments,
instead of just $j(n)$ in the general case; see Theorem~\ref{thm2} below.
We remark that for the special case of real symmetric second-order systems
and expansion point $s_0 = 0$, this result can be traced back to~\cite{SuC91}.

In this section, we first give an exact definition of Hermitian 
special second-order systems and higher-order systems, and then we
prove the stronger moment-matching property stated in Theorem~\ref{thm2}.

\subsection{Hermitian special second-order systems}

We say that a special second-order system~\mr{ASA} is \textit{Hermitian}
if the matrices in~\mr{ASA} and~\mr{AF1}, respectively~\mr{AF2},
satisfy the following properties:
\beq{HAC}{
L=B^H,\q P_0 = P_0^H,\q P_1 = P_1^H,\q F_1 = F_2,\q G=G^H.
}
Recall that RCL circuits are described by special second-order 
systems of the form~\mr{ACC} with real matrices defined in~\mr{ACD}.
Clearly, these systems are Hermitian.

Using~\mr{HAA},~\mr{HAC}, and~\mr{AF1}, respectively~\mr{AF2}, one readily verifies that 
the data matrices~\mr{CAM}, respectively~\mr{CAO}, of the equivalent first-order 
formulation~\mr{AFA2} satisfy the relations
\beq{HAE}{
\bea{rl}
\Jc\, \bigl(s_0\, \Ec - \Ac \bigr) & = \bigl(s_0\, \Ec - \Ac \bigr)^H \Jc,\q 
\Jc\, \Ec = \Ec\, \Jc,\q \Jc = \Jc^H,\\[6pt]
\Lc^H & = \Jc\, \Bc,
\ea
}
where
\[
\Jc := \Matrix{I_N & 0 \cr
              0 & -I_{N_0}}.  
\] 
Since the reduced-order model is structure-preserving, the data matrices~\mr{XXX1}
satisfy analogous relations.
More precisely, we have
\beq{HAX}{
\bea{rl}
\Jc_n\, \bigl(s_0\, \Ec_n - \Ac_n \bigr) 
   & = \bigl(s_0\, \Ec_n - \Ac_n \bigr)^H \Jc_n,\q 
       \Jc_n\, \Ec_n = \Ec_n\, \Jc_n,\q \Jc_n = \Jc_n^H,\\[6pt]
\Lc_n^H & = \Jc_n\, \Bc_n,
\ea
}
where
\[
\Jc_n : = \Matrix{I_n & 0 \cr
              0 & -I_n}.   
\] 

\subsection{Hermitian higher-order systems}

We say that a higher-order system~\mr{AHA} is \textit{Hermitian}
if the matrices in~\mr{AHA} satisfy the following properties:
\beq{HBA}{
P_i = P_i^H,\q 0\leq i\leq l,\q L_0 = B^H,\q
L_j = 0,\q 1\leq j \leq l-1.
}
  
In this case, we define matrices
\[
\hat{P}_j := \sum_{i=0}^{l-j} s_0^i P_{j+i},\q j = 0,1,\dots,l,
\]
and set
\beq{JAK2}{
\Jc :=  \Matrix{  I & -s_0 I & 0 & \cdots & 0\cr
                  0 &      I & -s_0 I & \ddots & \vdots \cr
             \vdots & \ddots &  \ddots  & \ddots  & 0 \cr
                  0 & \cdots & 0 & I & -s_0 I\cr
                       0 & 0 & \cdots & 0 & I }
 \Matrix{\hat{P}_1 & \hat{P}_2 &  \cdots & \hat{P}_{l-1} & I \cr
                \hat{P}_2 & \Addots & \Addots &  \hat{P}_l & 0 \cr
                 \vdots & \Addots & \Addots &  0 & \vdots \cr
                \hat{P}_{l-1} & \Addots & \Addots & \vdots & \vdots \cr
                \hat{P}_l & 0 & \cdots & 0 & 0}.
}
Note that, in view of~\mr{HBA}, we have
\beq{HBB}{
\hat{P}_j = \hat{P}_j^H,\q j=0,1,\dots,l.
}
Using~\mr{HBA}--\mr{HBB}, one can verify that the data matrices $\Ac$, $\Ec$,
$\Bc$, $\Lc$ defined in~\mr{CAC} satisfy the following relations:
\beq{HAM}{
\Jc\, \bigl(s_0\, \Ec - \Ac \bigr) = \bigl(s_0\, \Ec - \Ac \bigr)^H \Jc,\q 
\Jc\, \Ec = \Ec^H \Jc,\q \Lc^H = \Jc\, \Bc.
}
Since the reduced-order model is structure-preserving, the data matrices~\mr{XXX1}
satisfy the same relations.
More precisely, we have
\beq{HAY}{
\bea{rl}
\Jc_n\, \bigl(s_0\, \Ec_n - \Ac_n \bigr) & = \bigl(s_0\, \Ec_n - \Ac_n \bigr)^H \Jc_n,\q 
\Jc_n\, \Ec_n = \Ec_n^H \Jc_n,\\[6pt]
\Lc_n^H & = \Jc_n\, \Bc_n,
\ea
}
where $\Jc_n$ is defined in analogy to $\Jc$.
 
\subsection{Key relations}

Our proof of the enhanced moment-matching property in the Hermitian case
is based on some key relations that hold true for both special second-order
and higher-order systems.
In this subsection, we state these key relations.

Recall the definition of the matrix $\Mc$ in~\mr{PTR2}.
The relations~\mr{HAE}, respectively~\mr{HAM}, readily imply the following identity:
\beq{HAG}{
\Mc^H \Jc = \Jc\, \Ec\, \bigl(s_0\, \Ec - \Ac \bigr)^{-1}.
}
It follows from~\mr{HAG} that
\beq{HAJ}{
\bigr(\Mc^H\bigl)^i\, \Jc =
 \Jc\, \Bigl(\Ec\, \bigl( s_0\, \Ec - \Ac \bigr)^{-1}\Bigr)^i,\q
i=0,1,\dots\, .
}
Similarly, the relations~\mr{HAX}, respectively~\mr{HAY}, imply 
\[
\Mc_n^H \Jc_n = \Jc_n\, \Ec_n\, \bigl(s_0\, \Ec_n - \Ac_n \bigr)^{-1}.
\]
It follows that
\beq{HAN}{
\bigr(\Mc_n^H\bigl)^i\, \Jc =
 \Jc_n\, \Bigl(\Ec_n\, \bigl( s_0\, \Ec_n - \Ac_n \bigr)^{-1}\Bigr)^i,\q
i=0,1,\dots\, .
}
Also, recall from~\mr{HAE}, respectively~\mr{HAM}, that
\beq{HHH1}{
\Lc^H = \Jc\, \Bc,
}
and from~\mr{HAX}, respectively~\mr{HAY}, that
\beq{HHH2}{
\Lc_n^H = \Jc_n\, \Bc_n.
}
Finally, one readily verifies the following relation:
\beq{SPX18}{
\Vc_n^H\, \Jc\, \Ec\, \Vc_n  = \Jc_n\, \Ec_n.
}

\subsection{Matching twice as many moments}

In this subsection, we present our enhanced version of Theorem~\ref{thm1}
for the case of Hermitian special second-order or higher-order systems.

First, we establish the following proposition.

\begin{proposition} \label{prop1}
Let\/~$n = n(j)$ be of the form\/~\mr{KKE}.
Then, the data matrices~\mr{XXX1} of the structure-preserving Pad\'e-type
model satisfy
\beq{SPX10}{
\Lc\, \Mc^i\, \Vc_n = \Lc_n\, M_n^i\q \mbox{for all}\q i=0,1,\dots,j.
}
\end{proposition}

\begin{proof}
Recall that $\Lc_n = \Lc\, \Vc_n$.
Thus~\mr{SPX10} holds true for $i=0$.

Let $1\leq i \leq j$. 
In view of~\mr{HAJ}, we have
\[
\bigr(\Mc^H\bigl)^i\, \Jc = 
 \Jc\, \Bigl(\Ec\, \bigl( s_0\, \Ec - \Ac \bigr)^{-1}\Bigr)^i.
\]
Together with~\mr{HHH1}, it follows that
\[
\bigr(\Mc^H\bigl)^i\, \Lc^H = \bigr(\Mc^H\bigl)^i\, \Jc\, \Bc =
 \Jc\, \Bigl(\Ec\, \bigl( s_0\, \Ec - \Ac \bigr)^{-1}\Bigr)^i \Bc.
\]
Since $\bigl( s_0\, \Ec - \Ac \bigr)^{-1} \Bc = \Rc$, it follows
that
\[
\bigr(\Mc^H\bigl)^i\, \Lc^H = \Jc\, \Ec\, 
 \Bigl(\bigl( s_0\, \Ec - \Ac \bigr)^{-1} \Ec \Bigr)^{i-1} \Rc
= \Jc\, \Ec\, \Mc^{i-1} \Rc.
\]
Using~\mr{PTR10} (with $i$ replaced by $i-1$),~\mr{SPX18},~\mr{HAN}, and~\mr{HHH2},
we obtain 
\[
\bea{rl}
\Vc_n^H\, \bigr(\Mc^H\bigl)^i \Lc^H  
  & = \Vc_n^H\, \Jc\, \Ec\, \bigl(\Mc^{i-1} \Rc \bigr) \\[4pt]
  & = \Vc_n^H\, \Jc\, \Ec\, \Vc_n\, \Mc_n^{i-1} \Rc_n \\[4pt]
  & = \bigl(\Vc_n^H\, \Jc\, \Ec\, \Vc_n\bigr) \bigl(\Mc_n^{i-1} \Rc_n\bigr) \\[4pt]
  & = \Jc_n\, \Ec_n\, \Mc_n^{i-1} \Rc_n \\[4pt]
  & = \Jc_n\, \Ec_n\, \Mc_n^{i-1} \bigl( s_0\, \Ec - \Ac \bigr)^{-1} \Bc_n \\[4pt]
  & = \Jc_n\, \Bigl(\Ec_n\, \bigl( s_0\, \Ec - \Ac \bigr)^{-1}\Bigr)^i \Bc_n \\[4pt]
  & = \bigl(\Mc_n^H\bigr)^i \Jc_n\, \Bc_n = \bigl(\Mc_n^H\bigr)^i \Lc_n^H.
\ea
\]
Thus the proof is complete.
\end{proof} 

The following theorem contains the main result of this section.

\begin{theorem} \label{thm2}
Let\/~$n = n(j)$ be of the form\/~\mr{KKE}.
In the Hermitian case, the structure-preserving Pad\'e-type
model defined by the data matrices~\mr{XXX1} satisfies:
\[
H_n(s) = H(s) + {\Oc}\bigl((s-s_0)^{2 j(n)} \bigr).
\]
\end{theorem}

\begin{proof}
Let $j=j(n)$.
We need to show that
\beq{HCC}{
\Lc\, \Mc^i\, \Rc = \:c_n\, \Mc_n^i\, \Rc_n\q \mbox{for all}\q i=0,1,\dots,2j-1.
}
By~\mr{PTR10} and~\mr{SPX10}, we have
\[
\bea{rl}
\Lc\, \Mc^{i_1 + i_2}\, \Rc 
   & =  \bigl(\Lc\, \Mc^{i_1}\bigr) \bigl( \Mc^{i_2}\, \Rc \bigr)\\[4pt]
   & =  \bigl(\Lc\, \Mc^{i_1}\bigr) \bigl(\Vc_n\, \Mc_n^{i_2}\, \Rc_n \bigr)\\[4pt]  
   & =  \bigl(\Lc\, \Mc^{i_1}\, \Vc_n\bigr) \bigl(\Mc_n^{i_2}\, \Rc_n \bigr)\\[4pt]  
   & =  \bigl(\Lc_n\, \Mc_n^{i_1}\bigr) \bigl(\Mc_n^{i_2}\, \Rc_n \bigr)
     =  \Lc_n\, \Mc_n^{i_1+i_2}\, \Rc_n 
\ea
\] 
for all $i_1=0,1,\dots,j-1$ and $i_2=0,1,\dots,j$.
This is just the desired relation~\mr{HCC}, and thus the proof is complete.
\end{proof}

\section{The SPRIM algorithm} \label{sec-sprim}

In this section, we apply the machinery of structure-preserving Pad\'e-type
reduced-order modeling to the class of Hermitian special second-order systems that
describe RCL circuits.

Recall from Section~\ref{sec-2nd} that a first-order formulation of  
RCL circuit equations is given by~\mr{ABC} with data matrices defined in~\mr{ABD}.
Here, we consider first-order systems~\mr{ABC} with data matrices of the
slightly more general form
\beq{SSS1}{
\Ac = \Matrix{ -P_0 & -F \cr F^H & 0},\q
\Ec = \Matrix{ P_1 & 0 \cr 0 & G},\q
\Bc = \Matrix{ B \cr 0}.
}
Here, it is assumed that the subblocks $P_0$, $P_1$, and $B$ have the
same number of rows, and that the subblocks of $\Ac$ and $\Ec$
satisfy $P_0 \succeq 0$, $P_1  \succeq 0$, and $G \succ 0$.
Note that systems~\mr{ABC} with matrices~\mr{SSS1} are in particular
Hermitian.
Furthermore, the transfer function of such systems is given by
\[
H(s) = \Bc^H \bigl(s\, \Ec - \Ac\bigr)^{-1} \Bc.
\]

The PRIMA algorithm~\cite{OdaCP97,OdaCP98} is a reduction technique
for first-order systems~\mr{ABC} with matrices of the form~\mr{SSS1}.
PRIMA is a projection method that uses suitable basis matrices for the
block-Krylov subspaces $\Kc_n\bigl(\Mc,\Rc\bigr)$; see~\cite{Fre99d}.
More precisely, let $\hat{\Vc}_n$ be any matrix 
satisfying~\mr{QQQ1} and \mr{QQQ2}.
The corresponding $n$-th PRIMA model is then given by the projected
data matrices
\[
\hat{\Ac}_n := \hat{\Vc}_n^H \hat{\Ac}\,  \hat{\Vc}_n,\q
\hat{\Ec}_n := \hat{\Vc}_n^H \hat{\Ec}\,  \hat{\Vc}_n,\q
\hat{\Bc}_n := \hat{\Vc}_n^H \hat{\Bc}.
\]
The associated transfer function is
\[
\hat{H}_n(s) = \hat{\Bc}_n^H \bigl(s\, \hat{\Ec}_n - \hat{\Ac}_n\bigr)^{-1} 
                 \hat{\Bc}_n.
\]
For $n$ of the form~\mr{KKE}, the PRIMA transfer function satisfies
\beq{RED1}{
\hat{H}(s) = H(s) + {\cal O}\left( (s-s_0)^{j(n)}\right).
}

Recently, we introduced the SPRIM algorithm~\cite{Fre04c} as a structure-preserving
and more accurate version of PRIMA.
SPRIM employs the matrix $\Vc_n$ obtained from $\hat{\Vc}_n$ via the 
construction~\mr{QQQ3} and~\mr{QQQ4}.
The corresponding $n$-th SPRIM model is then given by the projected
data matrices
\[
\Ac_n := \Vc_n^H \Ac\,  \Vc_n,\q
\Ec_n := \Vc_n^H \Ec\,  \Vc_n,\q
\Bc_n := \Vc_n^H \Bc.
\] 
The associated transfer function is
\[
H_n(s) = \Bc_n^H \bigl(s\, \Ec_n - \Ac_n\bigr)^{-1} 
                 \Bc_n.
\]
In view of Theorem~\ref{thm2}, we have
\[
H(s) = H(s) + {\cal O}\left( (s-s_0)^{2j(n)}\right),
\]
which suggests that SPRIM is ``twice'' as accurate as PRIMA. 

An outline of the SPRIM algorithm is as follows.

\begin{algorithm} \textup{(SPRIM algorithm for special second-order systems)}
\begin{itemize}
\item
Input: matrices
\[
\Ac = \Matrix{ -P_0 & -F \cr F^H & 0},\q
\Ec = \Matrix{ P_1 & 0 \cr 0 & G},\q
\Bc = \Matrix{ B \cr 0},
\]
where the subblocks $P_0$, $P_1$, and $B$ have the
same number of rows, and the subblocks of $\Ac$ and $\Ec$
satisfy $P_0 \succeq 0$, $P_1  \succeq 0$, and
$G \succ 0$;\\
an expansion point $s_0\in \real$.
\item Formally set
\[
\Mc = \left(s_0\, \Ec - \Ac \right)^{-1} \Cc, \q
\Rc = \left(s_0\, \Ec - \Ac\right)^{-1} \Bc.
\]
\item Until $n$ is large enough, run your favorite block-Krylov
subspace method $($applied to $\Mc$ and $\Rc$$)$ to construct the columns
of the basis matrix
\[
\hat{\Vc}_n = \Matrix{v_1 & v_2 & \cdots & v_n}
\]
of the $n$-th block-Krylov
subspace $\Kc_n\bigl(\Mc,\Rc\bigr)$, i.e.,
\[
\Span \hat{\Vc}_n = \Kc_n\bigl(\Mc,\Rc\bigr).
\]
\item Let
\[
\hat{\Vc}_n = \Matrix{V_1 \cr V_2}
\]
be the partitioning of $\hat{\Vc}_n$ corresponding to the block sizes
of $\Ac$ and $\Ec$.
\item Set
\[
\bea{rl}
\tilde{P}_0 & = V_1^H P_1 V_1,\q
\tilde{F} = V_1^H F V_2,\q
\tilde{P}_1 = V_1^{H} P_1 V_1,\q
\tilde{G} = V_2^H G V_2,\\[6pt]
\tilde{B} & = V_1^H B,
\ea
\]
and 
\[
\Ac_n = \Matrix{ -\tilde{P}_0 & -\tilde{F} \cr \tilde{F}^H & 0},\q
\Ec_n = \Matrix{ \tilde{P}_1 & 0 \cr 0 & \tilde{G}},\q
\Bc_n = \Matrix{ \tilde{B} \cr 0}.
\]
\item Output: the reduced-order model $\widetilde{H}_n$
in first-order form
\beq{RED32}{
H_n(s) = \Bc_n^H \bigl(s\, \Ec_n - \Ac_n\bigr)^{-1} 
                 \Bc_n
}
and in second-order form 
\beq{sprim_2nd_a}{
H_n(s) =\tilde{B}^H \left(s\, \tilde{P}_1 +
    \tilde{P}_0
+ \frac{1}{s} \tilde{F} \tilde{G}^{-1} \tilde{F}^H\right)^{-1} \tilde{B}.
} 
\end{itemize}
\end{algorithm} 

We remark that the main computational costs of the SPRIM
algorithm is running the block Krylov subspace method
to obtain $\hat{\Vc}_n$.
This is the same as for PRIMA.
Thus generating the PRIMA reduced-order model $\hat{H}_n$ and the
SPRIM reduced-order model $H_n$ involves
the same computational costs. 
 
On the other hand, when written in first-order form~\mr{RED32},
it would appear that the SPRIM model has state-space dimension $2n$, and thus
it would be twice as large as the corresponding PRIMA model.
However, unlike the PRIMA model, the SPRIM model can always
be represented in special second-order form~\mr{sprim_2nd_a};
see Subsection \ref{ssec-pres}.
In~\mr{sprim_2nd_a}, the matrices $\tilde{P}_1$, $\tilde{P}_0$,
and $ \tilde{P}_{-1} := \tilde{F} \tilde{G}^{-1} \tilde{F}^H$
are all of size $n\times n$, and the matrix $\tilde{B}$ is
of size $n\times m$.
These are the same dimensions as in the PRIMA model~\mr{RED1}.
Therefore, the SPRIM model $H_n$
(written in second-order form~\mr{sprim_2nd_a}) and of the corresponding
PRIMA model $\hat{H}_n$ indeed have the same state-space dimension~$n$.

\section{Numerical examples} \label{sec-examples}

In this section, we present results of some numerical experiments
with the SPRIM algorithm for special second-order systems.
These results illustrate the higher accuracy of the SPRIM reduced-order
models vs. the PRIMA reduced-order models.
 
\subsection{A PEEC circuit} 
 
\begin{figure}
\centering
\epsfig{file=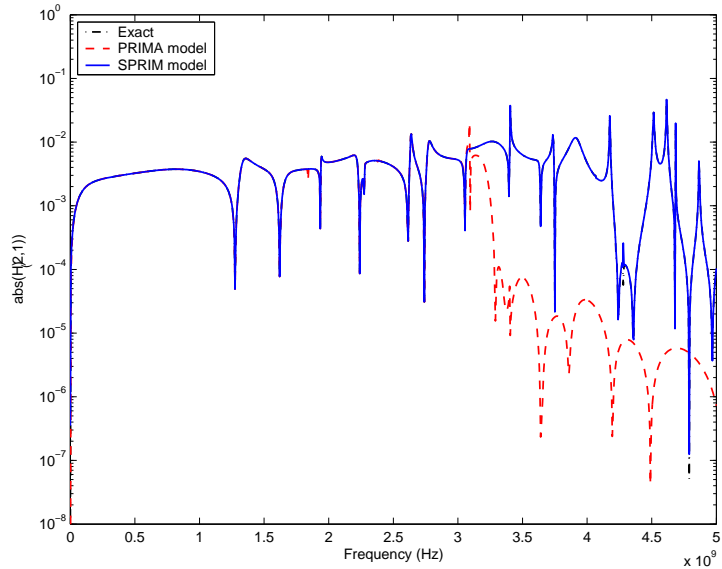,height=3in}
\caption{$\left|H_{2,1}\right|$ for PEEC circuit}
\label{peec}
\end{figure}

The first example is a circuit resulting from the so-called PEEC
discretization~\cite{Rue74} of an electromagnetic problem.
The circuit is an RCL network consisting of 2100 capacitors,
172 inductors, 6990 inductive couplings, and a single resistive source
that drives the circuit.
The circuit is formulated as a 2-port.
We compare the PRIMA and SPRIM models corresponding to the same
dimension $n$ of the underlying block Krylov subspace.
The expansion point $s_0 = 2\pi \times 10^{9}$ was used.
In Figure~\ref{peec}, we plot the absolute value of the $(2,1)$
component of the $2\times 2$-matrix-valued transfer function over
the frequency range of interest.
The dimension $n=120$ was sufficient for SPRIM to match
the exact transfer function.
The corresponding PRIMA model of the same dimension, however,
has not yet converged to the exact transfer function in large parts
of the frequency range of interest.
Figure~\ref{peec} clearly illustrates the better approximation
properties of SPRIM due to matching of twice as many moments
as PRIMA.

\subsection{A package model}
 
\begin{figure}
\centering
\epsfig{file=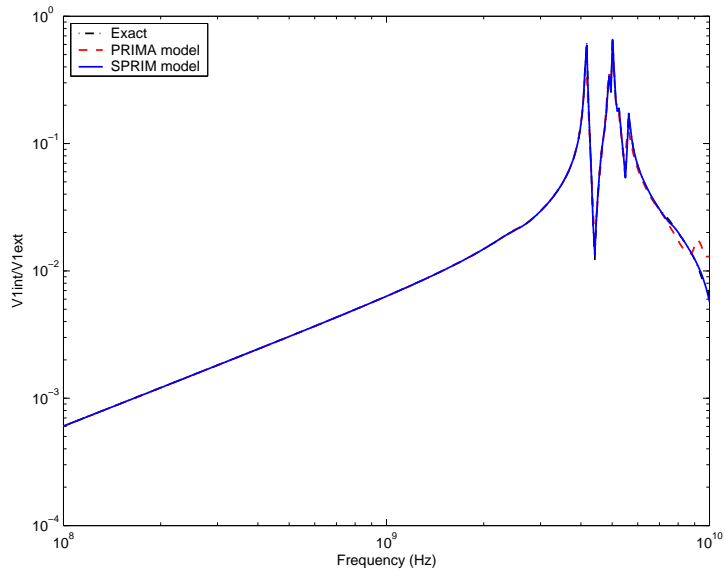,height=3in}
\caption{The package model}
\label{pack}
\end{figure}

\begin{figure}
\centering
\epsfig{file=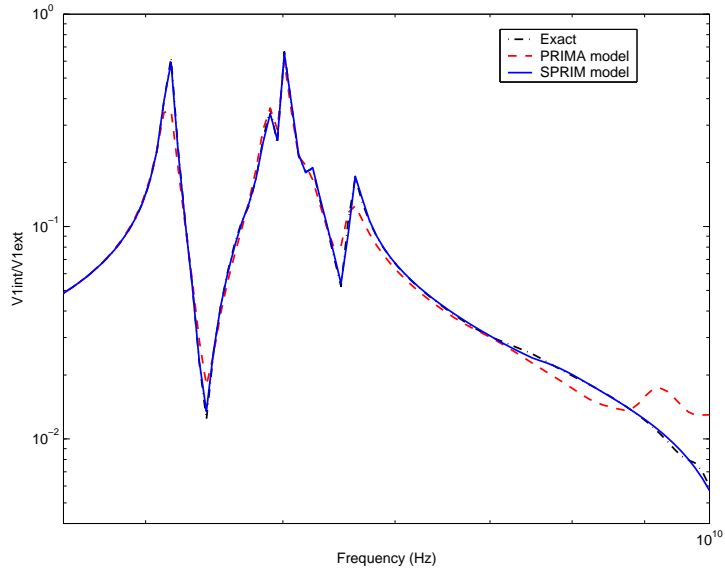,height=3in}
\caption{The package model, high frequencies}
\label{pack_det}
\end{figure}

The second example is a 64-pin package model used
for an RF integrated circuit. Only eight of the package pins carry
signals, the rest being either unused or carrying supply voltages.
The package is characterized as a 16-port component
(8 exterior and 8 interior terminals).
The package model is described by approximately 4000 circuit
elements, resistors, capacitors, inductors, and inductive
couplings.
We again compare the PRIMA and SPRIM models corresponding to the same
dimension $n$ of the underlying block Krylov subspace.
The expansion point $s_0 = 5\pi \times 10^{9}$ was used.
In Figure~\ref{pack}, we plot the absolute value of one of
the components of the $16\times 16$-matrix-valued transfer
function over
the frequency range of interest.
The state-space dimension $n=80$ was sufficient for SPRIM to match
the exact transfer function.
The corresponding PRIMA model of the same dimension, however,
does not match the exact transfer function very well near
the high frequencies; see Figure~\ref{pack_det}.
 
\subsection{A mechanical system}
 
\begin{figure}
\centering
\epsfig{file=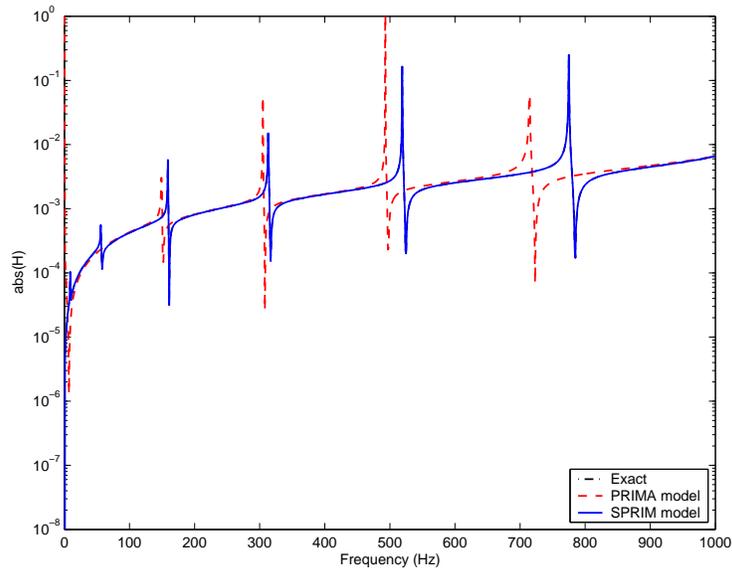,height=3in}
\caption{A mechanical system}
\label{shaft}
\end{figure}

Exploiting the equivalence (see, e.g.,~\cite{LozBEM00}) between
RCL circuits and mechanical systems, both PRIMA and SPRIM can
also be applied to reduced-order modeling of mechanical systems.
Such systems arise for example in the modeling and simulation
of MEMS devices.
In Figure~\ref{shaft}, we show a comparison of PRIMA and SPRIM
for a finite-element model of a shaft.
The expansion point $s_0 = \pi \times 10^{3}$ was used.
The dimension $n=15$ was sufficient for SPRIM to match
the exact transfer function in the frequency range of interest.
The corresponding PRIMA model of the same dimension, however,
has not converged to the exact transfer function in large parts
of the frequency range of interest.
Figure~\ref{shaft} again illustrates the better approximation
properties of SPRIM due to the matching of twice as many moments
as PRIMA.

\section{Concluding remarks} \label{sec-cremarks}
  
We have presented a framework for constructing structure-preserving
Pad\'e-type reduced-order models of higher-order linear dynamical
systems.
The approach employs projection techniques and Krylov-subspace machinery for
equivalent first-order formulations of the higher-order systems.
We have shown that in the important case of Hermitian higher-order
systems, our structure-preserving Pad\'e-type model reduction is twice
as accurate as in the general case.
Despite this higher accuracy, the models produced by our approach are
still not optimal in the Pad\'e sense. 
This can be seen easily by comparing the degrees of freedom of general
higher-order reduced models of prescribed state-space dimension, with
the number of moments matched by the Pad\'e-type models generated by
our approach.
Therefore, structure-preserving true Pad\'e model reduction remains
an open problem.

Our approach generates reduced models in higher-order form via 
equivalent first-order formulations.
It would be desirable to have algorithms that construct the same
reduced-order models in a more direct fashion, without the detour
via first-order formulations.
Another open problem is the ``optimal'' way of constructing
basis vectors for the structured Krylov subspaces that arise for the 
equivalent first-order formulations.
In particular, an algorithm for this task should be both
computationally efficient and numerically stable.
Some related work on this problem is described in the recent
report \cite{Li04}, but many questions remain open.
Finally, the proposed approach is a projection technique,
and as such, it requires the storage of all the vectors
used in the projection.
This clearly becomes an issue for systems with very large
state-space dimension.

\bibliographystyle{abbrv}  

\bibliography{../siam_book/biblio}
\end{document}